\documentclass[12pt]{article}

\usepackage[T2A]{fontenc}
\usepackage[utf8]{inputenc}

\usepackage{amsfonts}
\usepackage{amssymb}
\usepackage{amsmath}
\usepackage{amsthm}

\def \le {\leqslant}
\def \ge {\geqslant}

\begin{document}

\begin{Large}
\centerline{\bf Badly approximable vectors in affine subspaces: }
\centerline{\bf Jarn\'{\i}k-type result }
\end{Large}
 \vskip+0.5cm \centerline{by {\bf Nikolay Moshchevitin}.\footnote{supported by the grant RFBR № 09-01-00371a
 }  }

 \vskip+0.5cm
\begin{small}
{\bf Abstract.}\,
Consider irrational affine subspace $ A\subset \mathbb{R}^d$ of dimension $a$.
We prove that the set
$$
\{\xi =(\xi_1,...,\xi_d)  \in { A}:\,\,\,
\
q^{1/a}\cdot
\max_{1\le i \le d} ||q\xi_i||  \to \infty,\,\,\,\,
q\to \infty\}
$$
is an $\alpha$-winning set for every $\alpha \in (0,1/2]$

\end{small}

\vskip+0.5cm

This simple short communication may be considered as a supplement to our short paper \cite{m}.

{\bf 1.\, Jarn\'{\i}k's result in simultaneous Diophantine approximations.}
\,\,\,
All numbers in this paper are real. Notation $||\cdot||$ stands for the distance to the nearest integer.
In 1938
  V. Jarn\'{\i}k (see \cite{TB}, Satz 9 and \cite{J}, Statement {\bf E}) proved the following result.

 {\bf Theorem 1.} (V. Jarn\'{\i}k)\,\,\,
{\it Suppose that among numbers $\xi_1,...,\xi_d$ there are at least two numbers which are linearly
independent over $\mathbb{Z}$, together with 1.
Then
$$
\limsup_{t\to +\infty}
\,\,\left(t \cdot
\min_{q\in \mathbb{Z},\, 1\le q \le t}
\,\,
\max_{1\le i \le d} \,\,
||q\xi_i ||\right) = +\infty
.
$$
}

Nothing more can be said in a general situation.
In 1926
A. Khintchine \cite{HINS} proved the following result.

{\bf Theorem 2.} (A. Khintchine) \,\,\,{\it
Let
 $\psi(t)$
increases to  infinity as  $t\to +\infty$.
Then there exist two algebraically independent
real numbers $\xi_1,\xi_2$ such that
for all $t$ large enough one has
$$
t \cdot
\min_{ q\in \mathbb{Z},\, 1\le q \le t} \, \max_{i=1,2}||q\xi_i ||\le \psi (t).
 $$}

A general form of such a result one can find in
Jarn\'{\i}k's paper \cite{J}. A corresponding $\limsup$ result is due to J. Lesca \cite{LL}:

{\bf Theorem 3.} (J. Lesca)\,\,\,
{\it Let $d\ge 2$.  Let $\psi(t)$ be a positive continuous  function in $t$ such that the function
$t\mapsto \psi (t) /t$ is a decreasing function.
 Suppose that
$$
\limsup_{t\to\infty}{\psi(t)}=+\infty.
$$
Then the set of all vectors $\xi = (\xi_1,...,\xi_d)\in \mathbb{R}^d$,
 containing of algebraically independent elements,
 such that
$$
 t \cdot
\min_{ q\in \mathbb{Z},\, 1\le q \le t} \, \max_{1\le i\le d}||q\xi_i ||\le \psi (t)
$$
for all $t$ large enough,
  being intersected with a given open set
  ${\cal G} \subset \mathbb{R}^{d}$ is of cardinality continuum.
}

We would like to note that Jarn\'{\i}k's Theorem 1 as well as some
other theorems by Khintchine and  V. Jarn\'{\i}k were discussed
and generalized in author's survey \cite{MM}. In particular in
\cite{MM}, Section 4.1 (see also \cite{MMM}) one can find an
improvement of Theorem 1 in terms of the best approximation
vectors.

{\bf 2.\, Affine subspaces.}
\,\,\,
Let
$\mathbb{R}^d$
be a Euclidean space with the coordinates $(x_1,...,x_d)$,
let $\mathbb{R}^{d+1}$
be a Euclidean space with the coordinates $(x_0,x_1,...,x_d)$.
 Consider an affine subspace $A\subset \mathbb{R}^d$.
Let $a  = {\rm dim }\, A \ge 1$. Define
the affine subspace
${\cal A}\subset \mathbb{R}^{d+1}$ in the following way:
$$
{\cal A} =\{ {\bf x} = (1,x_1,...,x_d):\,\, (x_1,...,x_d) \in A\}.
$$
 We define {\it linear} subspace
$
\hbox{\got A} = {\rm span} \,{\cal A},
$
as the smallest linear subspace in $\mathbb{R}^{d+1}$ containing
${\cal A}$.
So ${\rm dim}\,  \hbox{\got A} = a+1$.

Consider a sublattice $ \Gamma ( A) = \hbox{\got A} \cap \mathbb{Z}^{d+1}$
of the integer lattice $\mathbb{Z}^{d+1}$. We see that
$$
0\le {\rm dim}\,  \Gamma ( A)
\le a+1.
$$
Of course here for a lattice $\Gamma \subset \mathbb{Z}^{d+1}$ by ${\rm dim}\,  \Gamma$
we mean the dimension of the linear subspace
${\rm span}\, \Gamma$.

In the case $ {\rm dim}\,  \Gamma ( A)
= a+1 ={\rm dim}\,  \hbox{\got A}$ we define  $ A$
to be a {\it completely rational} affine subspace in $\mathbb{R}^{d}$.
For a completely rational affine subspace $A$ by $d(A)$ we denote the fundamental $(a+1)$-dimensional volume of
the lattice $\Gamma (A)$.


We see from Dirichlet principle that for any completely rational
affine subspace $ A$ of dimension $a$ there exists a positive
constant $\gamma = \gamma ({ A})$ such that for any $\xi
=(\xi_1,...,\xi_d) \in A$
 the inequality
$$
 \max_{1\le i \le d} ||q\xi_i|| \le \frac{\gamma}{ q^{1/a}}
$$
has infinitely many solutions in positive integers $q$.

One can easily see that for any affine subspace ${ A}$ of dimension $a$
the set
$$\Omega =
\{\xi =(\xi_1,...,\xi_d)  \in { A}:\,\,\,
\inf_{q\in \mathbb{Z}_+ }\, q^{1/a} \cdot \max_{1\le i \le d} ||q\xi_i ||\,\,\, >0\}
$$
is an 1/2-winning set in $ A $. Here we do not want to discuss the
definitions $(\alpha, \beta)$-games and $(\alpha,\beta)$-winning
set or $\alpha$-winning set. This definitions were given in W.M.
Schmidt's paper \cite{sss}. All  the definitions and basic
properties of wining sets one can find in  the book \cite{SC},
Chapter 3. In particular, every $\alpha$-winning set in $A$ has
full Hausdorff dimension. A countable intersection of
$\alpha$-winning sets inn $A$ is also an $\alpha$-winning set.

In the case when ${ A}$ is not a completely rational subspace the
result about winning property of the set $\Omega$ admits a small improvement. This improvement is related to
Jarn\'{\i}k's result cited behind.

{\bf Theorem 4.}\,\,{\it
Let
$$
0< \alpha < 1,\,\,\, 0< \beta < 1,\,\,\, \gamma= 1+\alpha \beta - 2\alpha >0 .$$
 Suppose that ${\rm dim}\,  \Gamma ( A)
< a$.
Then the set
\begin{equation}\label{mai}
\Omega^* =
\{\xi =(\xi_1,...,\xi_d)  \in { A}:\,\,\,
\
q^{1/a}\cdot
\max_{1\le i \le d} ||q\xi_i||  \to \infty,\,\,\,\,
q\to \infty\}
\end{equation}
is $(\alpha,\beta)$-winning set in $ A.
$
In particular, it is an $\alpha$-winning set for every $\alpha \in (0,1/2]$.

}

Here we should note that certain results concerning badly approximable vectors in affine subspaces
one can find in
\cite{K,Kl,m,z}.

{\bf 3. Lemmata.}\,\,\,\,
Consider the set of all $(a+1)$-dimensional complete sublattices of
the integer  lattice $\mathbb{Z}^{d+1}$. It is a countable set.
One can easily see that for any positive $H$  there exist not more than a finite number of such
sublattices $\Gamma$
with the fundamental volume ${\rm det} \, \Gamma \le H$.
Hence we can order the set
$\{ V_\nu\}_{\nu = 1}^\infty$
 of all $a$-dimensional  affine subspaces in $\mathbb{R}^d$ in such a way that
values
 $d_\nu =d(V_\nu) = {\rm det}\, \Gamma (V_\nu)$ form an increasing sequence:
$$
1 =
d_1
\le
d_2\le \cdots\le
d_\nu \le d_{\nu+1}\le \cdots
.$$
We see that
\begin{equation}\label{inf}
 { d}_\nu \to \infty,\,\,\, \nu \to \infty.
\end{equation}
Some of consecutive values of $d_\nu$ may be equal. We define a
sequence $d_{\nu_k}$ of all different elements from the sequence
$\{ d_\nu\}$:
$$
1 =
d_{\nu_1} = ... = d_{\nu_2 -1} < d_{\nu_2}= .... = d_{\nu_2-1}< d_{\nu_3 }
=
...
<
d_{\nu_k} =
 d_{\nu_k+1}
=
...
=
d_{\nu_{k+1}-1}
<
d_{\nu_k} = ...
$$
(of course $\nu_1 = 1$).
For $V_j$ we define the affine subspace
${\cal V}_j\subset \mathbb{R}^{d+1}$ as
$$
{\cal V}_j =\{ {\bf x} = (1,x_1,...,x_d):\,\, (x_1,...,x_d) \in V_j\}
$$
and consider  linear subspace
$
\hbox{\got V}_j = {\rm span} \,{\cal V}_j.
$

In the sequel  for $\xi = (x_1,...,x_d)\in A$ we consider $a$-dimensional ball
$$
B(\xi, \rho) = \{ \xi'=(\xi_1',...,\xi_d')\in A:\,\,\, \max_{1\le i\le d} |\xi_i - \xi'_i|\le \rho\}
$$
and $d$-dimensional ball
$$
\hat{B}(\xi, \rho) = \{ \xi'=(\xi_1',...,\xi_d')\in \mathbb{R}^d:\,\,\, \max_{1\le i\le d} |\xi_i - \xi'_i|\le \rho\}
.$$
Obviously
$$
B(\xi, \rho)  = \hat{B}(\xi, \rho)\cap A.
$$

{\bf Lemma 1.}\,\,{\it
Suppose that  $ U,  V\subset \mathbb{R}^{d}$ are two affine subspaces.
Put $  L =  U\cap  V$ and suppose that
$
 {\rm dim}\, U > {\rm dim}\,  L.
$
Suppose that affine subspace $   L' \subset   U$
has dimension
$ {\rm dim}\,  L' = {\rm dim}\,  U -1$, and $  L'\cap L = \varnothing$.
Define $\overline{U} \subset  U$ to be a half-subspace  with  the boundary
$ L'$ and such that $ \overline{ U}\cap  L = \varnothing$.
Then
$
{\rm dist} (\overline{ U},  V) >0.
$}

Proof.\,\,
In affine subspace $ {\rm aff}\, (  U\cup  V)$  of dimension
$ w =  {\rm dim}\,  U+ {\rm dim}\,  V - {\rm dim}\,  L$
there exists an affine subspace $  L'' \supset   L'$ with dimension
${\rm dim}\,  L'' = w-1$ such that $  L''\cap  V =\varnothing$.
So ${\rm dist} \, ( L'',  U)
>0$.
The subspace $ L''$ divides $ {\rm aff}\, (  U\cup  V)$  into two parts,
and lemma follows.$\Box$

{\bf Corollary.}\,\,{\it Consider two affine subspaces $A, V\subset \mathbb{R}^d$.
Suppose that for $ \xi \in A$
the ball $B(\xi,\rho)\subset A $
satisfies the property
$$
{\rm dist}
(B(\xi,\rho),A\cap V)  \ge \varepsilon >0 .
$$
Then there exists positive
 $\delta = \delta (A,V,\xi, \varepsilon)$ such that for any $ \xi ' \in {B(\xi,\rho)}$  one has
$$
\hat{B} (\xi', \delta) \cap V = \varnothing.
$$}

Proof. \,\, From the conditions of our Corollary we see that ${\rm
dim}\, (A\cap V) < {\rm dim }\, A.$ So we can take a subspace $L'$
of dimension $ {\rm dim}\, L' = {\rm dim}\, A -1$ which separates
the ball $B(\xi,\rho)$ from  the subspace $A\cap V$ in $A $. Now
we use
 Lemma 1.$\Box$

{\bf Lemma 2.}\,\,{\it
Let $\rho>0$ and $\xi \in A$.
Consider
a ball  $\hat{B}(\xi,\rho)\subset\mathbb{R}^d$
such that
\begin{equation}\label{dee}
 \hat{B}(\xi,\rho) \cap  \hbox{\got V}_j =\varnothing,\,\,\,\, 1\le j \le n.
\end{equation}
Define $k= k(n)$ from the condition
\begin{equation}\label{dee1}
 \nu_{k} \le n
< \nu_{k+1}
.
\end{equation}
Put
 \begin{equation}\label{si}
  \kappa =\kappa_{d,\xi} =
(2\sqrt{d})^a
\times \sqrt{1+ (|\xi_1|+1)^2 +...+ (|\xi_d|+1)^2}
 ,\,\,\,
\sigma =
\sigma_{a,d,\xi }=
\frac{1 }{\kappa_{d,\xi} (a+1)!}
\end{equation}
and
$$T = (\sigma d_{\nu_n} \rho^{-a})^{\frac{1}{a+1}}.
$$
Then the set of all rational points $\left(\frac{b_1}{q},...,\frac{b_d}{q}\right)\in \hat{B}(\xi,\rho)$
with $q\le   T$  lie in a certain $(a-1)$-dimensional affine subspace.}

Proof.\,\, We may suppose that
the set of rational points
from
$ \hat{B}(\xi,\rho)$
with $q\le T$
consists of more than $a$ points (otherwise there is nothing to prove).
We take arbitrary $a+1$ points
$$
\left(\frac{b_{1,j}}{q_j},...,\frac{b_{d,j}}{q_j}\right) \in \hat{B}(\xi,\rho),\,\,\,\,
1\le q_j \le T,\,\,\,\,
{\rm g.c.d.} (q_j,b_{1,j},...,b_{d,j}) =1,\,\,\,\,
1\le j \le a+1
$$
and  prove that primitive integer vectors
\begin{equation}\label{ve}
{\bf b}_j= (q_j,b_{1,j},...,b_{d,j}) ,\,\,\,\, 1\le j \le a+1
\end{equation}
are linearly dependent.
Then the lemma will be proved.

All integer vectors (\ref{ve})   belong to the cylinder
$$
C =C_{ \xi} (T,\rho) =
\{ {\bf x} = (x_0,x_1,....,x_d) \in \mathbb{R}^{d+1}:\,\,\,
0\le x_0 \le T,\,\,\, \max_{1\le j \le d} |x_0\xi_j -x_j|\le \rho T.
\}
$$
 Suppose that they are linearly independent.
Then $\hbox{\got L} = {\rm span} ({\bf b}_1,...,{\bf b}_{a+1})$
is an $(a+1)$-dimensional completely rational linear subspace.
By $D$ we denote the fundamental $(a+1)$-dimensional volume of the lattice $\hbox{\got L}\cap \mathbb{Z}^{d+1}$.
From (\ref{dee}) we see that
$$
\hbox{\got L} \neq \hbox{\got V}_j
,\,\,\,\, 1\le j \le n.
$$
From (\ref{dee1}) we see that
\begin{equation}\label{DE}
 D  \ge d_{\nu_n} .
\end{equation}
Now we consider the section
$
\hbox{\got L} \cap C
$
which is an $(a+1)$-dimensional convex polytope. As it is inside $C$, its $(a+1)$-dimensional measure is less than
$$
(2\sqrt{d}\rho T )^a
\times  T\sqrt{1+ (|\xi_1|+1)^2 +...+ (|\xi_d|+1)^2} =
\kappa\rho^a T^{a+1} =
\kappa  \sigma d_{\nu_n}.
$$
But the section $
\hbox{\got L} \cap C
$ consist of $a+1$ independent points from the lattice
 $\hbox{\got L}\cap \mathbb{Z}^{d+1}$. For the fundamental volume of this lattice we have lower bound (\ref{DE}).
That is why
$$
\frac{d_{\nu_n} }{(a+1)!} \le
\frac{D}{(a+1)!} <
\kappa \sigma d_{\nu_n}
=\frac{d_{\nu_n} }{(a+1)!}
.
$$
 This is a contradiction. Lemma is proved.$\Box$

 {\bf Lemma 3.}(W.M. Schmidt's  escaping lemma, Lemma 1B,  \cite{SC}, Chapter 3)  \,\,\,{\it
 Let $t$ be such that
 $$
 (\alpha \beta )^t <\frac{\gamma}{2}.
 $$
 Suppose   a ball $ B_j\subset A$ with the radius $\rho_j$ occurs in
the game (as a Black ball). Suppose $V$ is an $(d-1)$-dimensional affine subspace passing through the
 center of the ball $B_j$. Then White  can play in such a way that the ball $B_{j+t}$
 is contained in the halfspace $ \Pi$
 such that the boundary of
$ \Pi$ is parallel to the subspace $V$ and the distance between $\Pi$ and $V$ is equal to $\frac{\rho_j\gamma}{2}$.
 }

{\bf Corollary.} {\it  Suppose   a ball $ B_j\subset A$ with the radius $\rho_j$ occurs in
the game (as a Black ball).
Suppose that $  V, V' \subset A$ are two proper affine subspaces of $A$.
Then White  can play in such a way that the distance from the ball $B_{j+2t}$
to each of subspaces $V,V'$ is greater than $\frac{\rho_{j+t}\gamma}{2}$
(here $\rho_{j+t}$ is the radius of the ball $B_{j+t}$).
}

{\bf 4. Proof of Theorem 4.}
Suppose that $t = t(\alpha, \beta)$ satisfies the condition of Lemma 3. Put
$
j_k = 2tk$
and $R_0 = 1$.
Suppose that the first Black ball $B_0\subset A$ with the radius $\rho_0$ lies inside the box
$\{ \xi \in \mathbb{R}^d:\,\,\, \max_{1\le i \le d} |\xi_i| \le W\}$.
We shall
prove that
White
can play in such a way
that for any $\xi \in B_{j_r}$ one has
\begin{equation}\label{val}
\max_{1\le i \le d}|| q\xi_i || \ge \frac{(\alpha\beta)^t\gamma \rho_0}{2} \,  R^{-\frac{(a+1)r}{a}} \cdot q,\,\,\,\,
\forall
q < R_r
\end{equation}
with a certain $R_r$ which we define later in the inductive step.

We shall prove it by induction in $r$.

The base of induction is obvious.

Suppose that the ball $B_{j_{r-1}} = B(\xi_{j_{r-1}},\rho_{j_{r-1}})\in A,
\xi_{j_{r-1}} = (\xi_{j_{r-1},1},..., \xi_{j_{r-1},d})
$ which occurs as a Black ball
satisfies the condition specified.
Note that $ \rho_{j_{r-1}}= \rho_0 (\alpha\beta)^{j_{r-1}} $.
Consider the ball $\hat{B}_{j_{r-1}} = \hat{B}(\xi_{j_{r-1}},2\rho_{j_{r-1}})\in \mathbb{R}^d$.
Define $k_r$ as the maximal $k$ such that
$\hat{B}_{j_{r-1}} \cap \hbox{\got V}_j = \varnothing, 1\le j \le \nu_k$.
 Then we apply Lemma 2 to see that all  rational points
$\left(\frac{b_1}{q},...,\frac{b_d}{q}\right)\in \hat{B}_{j_{r-1}} $
with
$$
q\le
\left(
\sigma_{a,d,\xi_{j_{r-1}}} (2\rho_0)^{-a}\right)^{\frac{1}{a+1}}\left(\frac{1}{\alpha\beta}\right)^{\frac{2at(r-1)}{a+1}}
d_{\nu_{k_{r-1}}}^{\frac{1}{a+1}}
$$
 lie in a certain $(a-1)$-dimensional affine subspace. We denote this subspace  by $V_{r}'$.
As $\max_{1\le i \le d}|\xi_{j_{r-1},i}|\le W$ we see that
$$
\sigma_{a,d,\xi_{j_{r-1}}}
\ge \Sigma_{a,d,W} = \frac{1}{(2\sqrt{d})^a\,\sqrt{1+(W+1)^2d}\, (a+1)!}.
$$
We put
\begin{equation}\label{erer}
R_r
=
\left(
\Sigma_{a,d,W} (2\rho_0)^{-a}\right)^{\frac{1}{a+1}}\left(\frac{1}{\alpha\beta}\right)^{\frac{2at(r-1)}{a+1}}
d_{\nu_{k_{r-1}}}^{\frac{1}{a+1}}.
\end{equation}
By Corollary to Lemma 3 White can play in such a way that
\begin{equation}\label{11}
{\rm dist}
(B_{j_r},V_r)  \ge \frac{\gamma \rho_{j_{r-1}+t}}{2}
\end{equation}
and
\begin{equation}\label{22}{\rm dist}
(B_{j_r},V_r')  \ge \frac{\gamma \rho_{j_{r-1}+t}}{2}
\end{equation}
So the inductive step is described and we must show that
(\ref{val}) is valid. But it is clear from (\ref{22})  that for
any $\xi \in B_{j_r}$ one has
\begin{equation}\label{fle}
\max_{1\le i \le d}||q\xi_i || \ge
\frac{1}{2} \gamma \rho_{j_{r-1}+t} q  =
\frac{\gamma\rho_0}{2} \left({\alpha\beta}\right)^{(2t-1)r} q, \,\,\,\,
\forall
q< R_r.
\end{equation}
Moreover
by  Corollary to Lemma 1  from (\ref{11}) we see that
$$
k_r \to +\infty,\,\,\,\, r \to +\infty.
$$
  Hence
\begin{equation}\label{fle1}
d_{\nu_{k_r}} \to +\infty,\,\,\,\, r \to +\infty.
\end{equation}

Consider a point $\xi \in \bigcap_{j} B_j$.
For positive integer $
 q$
define $r$ from the condition
$$
R_{r-1} \le q < R_r.
$$
Then we make use of $\xi \in  B_{j_r}$.
From the inequality $q\ge R_{r-1}$ and (\ref{erer}) we see that
$$
\alpha \beta \ge \omega_1 q^{-\frac{a+1}{2atr}} d_{\nu_{k_{r-2}}}^{\frac{1}{2atr}},$$
where $ \omega_1 = \omega_1 (a,d,W,\alpha,\beta, t)>0$.
We substitute this estimate into
(\ref{fle}) to see that
$$
\max_{1\le i \le d}
||q\xi_i || \ge
\omega_2 q^{-1/a}
d_{\nu_{k_{r-2}}}^{1/a},\,\,\,\, R_{r-1} \le q < R_r,
 $$
with positive $ \omega_2 = \omega_2 (a,d,W,\alpha,\beta, t)$.
From (\ref{fle1}) for $\xi \in \bigcap_j B_j $
we deduce that
$$
q^{{1}/{a}}\cdot
\max_{1\le i \le d}||q\xi_i || \to +\infty,\,\,\,\, q \to \infty.
$$
 So White can enforce Black to reach
a point $\xi$ with the desired properties.
Theorem  4 is proved.$\Box$


\begin{thebibliography}{10}

\bibitem{TB}
V. Jarn\'{\i}k,\,\,\,
Zum Khintchineschen `` \"Ubertragungssatz'',
Travaux de L'Institut Mathematique de Tbilissi, T 3
(1938), 193 - 212.


\bibitem{J}
V. Jarn\'{\i}k,\,\,\,
 Eine Bemerkung uber diophantische approximationen, Math. Zeitschr. 72, 187 - 191 (1959).



\bibitem{HINS}
{A.Y. Khinchine},   \,\,\, \"Uber eine klasse linear Diophantine Approximationen, Rendiconti Circ.
Math. Palermo, 1926, 50, p.170 - 195.

\bibitem{K}
D. Kleinbock,\,\,
Extremal supspaces and their submanifolds, GAFA 13:2 (2003), 437 - 466.

\bibitem{Kl}
D. Kleinbock,\,\,
An extension of quantitative nondivergence and applications to Diophantine exponents,
Trans. Amer. Math. Soc. 360 (2008), 6497-6523.


\bibitem{LL}
J. Lesca,\,\,\,
 Sur un r\'esultat de Jarn\'\i k,  Acta Arith.
 11
(1966),  359--364.


\bibitem{MM}
N.G. Moshchevitin,\,\,\, Khintchine's singular Diophantine systems and their
applications, Russian Mathematical Surveys. 65:3 43 - 126 (2010).

\bibitem{MMM}
N.G. Moshchevitin, \,\,Best simultaneous approximations: Norms, signatures, and
asymptotic directions, Math. Notes 67:5 (2000), 618–624.


\bibitem{m}
N. Moshchevitin,\,\,
On Kleinbock's Diophantine result,
submited to Proc. Math. Debrecen, preprint available at arXiv:0906.1541 (2009, 2011).




\bibitem{sss}
 W.M. Schmidt,\,\,
 On badly approximable
         numbers and certain games,
         Trans. Amer. Math. Soc., 623 (1966), p. 178 -- 199.


\bibitem{SC}
 W.M. Schmidt,\,\, Diophantine Approximations, Lect. Not. Math.,  785 (1980).



\bibitem{z}
Y. Zhang,\,\,
Diophantine exponents of affine subspaces: The simultaneous approximation case,
J. Number Theory, 129:8 (2009),
 1976-1989.




\end{thebibliography}
\end{document}